\documentclass[11pt]{article}

\usepackage[T1]{fontenc}
\usepackage{lmodern}
\usepackage[a4paper,margin=31mm]{geometry}
\usepackage{amsmath,amssymb,amsthm,mathtools}
\usepackage{microtype}
\usepackage{enumitem}
\usepackage{xcolor}
\usepackage{tikz}
\usetikzlibrary{arrows.meta,calc,decorations.pathreplacing,patterns}
\usepackage[colorlinks=true,linkcolor=blue,citecolor=blue,urlcolor=blue]{hyperref}

\definecolor{revisionblue}{RGB}{0,92,170}

\numberwithin{equation}{section}

\newtheorem{theorem}{Theorem}

\newtheorem{lemma}[theorem]{Lemma}

\theoremstyle{remark}
\newtheorem{remark}[theorem]{Remark}
\newtheorem*{definition}{Definition}

\newcommand{\R}{\mathbb{R}}

\newcommand{\dist}{\operatorname{dist}}
\newcommand{\diam}{\operatorname{diam}}

\title{$p$-Dirichlet spaces over chord-arc domains}
\author{Huaying Wei\thanks{Center for Applied Mathematics and KL-AAGDM, Tianjin
University, Tianjin 300072, China.
\texttt{hywei@tju.edu.cn}}
\ and Michel Zinsmeister\thanks{Institut Denis Poisson, Universit\'e d'Orl\'eans, 
Orl\'eans, 45067, France.
\texttt{zins@univ-orleans.fr}}}
\date{}

\begin{document}

\maketitle
\insert\footins{\noindent\footnotesize
Research is supported  by the National Natural Science Foundation
of China (grant nos.~12271218 and 12571083).}

\begin{abstract}
Let $\Gamma$ be a rectifiable Jordan curve in the complex plane, let
$\Omega_i$ and $\Omega_e$ be its interior and exterior domains, respectively,
and let $1<p<\infty$.  Let $E$ be the vector space of restrictions to $\Gamma$
of functions in $C^1(\mathbb C)$.  We consider the following three seminorms on
$E$:
\begin{enumerate}
\item $\displaystyle \Vert u\Vert_i=\left(\frac{1}{2\pi}
\iint_{\Omega_i}|\nabla U_i(z)|^p\lambda_{\Omega_i}^{2-p}(z)\,dA(z)\right)^{1/p}$,
where $U_i$ is the harmonic extension of $u$ to $\Omega_i$ and
$\lambda_{\Omega_i}$ is the hyperbolic density of $\Omega_i$;
\item $\Vert u\Vert_e$, defined analogously on $\Omega_e$;
\item $\displaystyle \Vert u\Vert_{B_p(\Gamma)} =
\left(\frac{1}{4\pi^2}\iint_{\Gamma\times\Gamma}
\frac{|u(z)-u(\zeta)|^p}{|z-\zeta|^2}\,|dz|\,|d\zeta|\right)^{1/p}$.
\end{enumerate}
These three seminorms are known to be equivalent when $\Gamma$ is a
chord-arc curve.  We investigate the converse problem. 
\end{abstract}


\section{Introduction and results}
This paper continues our study \cite{WZ} of Dirichlet spaces over planar
domains.  Such spaces arise naturally in operator theory, planar potential
theory, and Teichm\"uller theory.  Here we pass from the quadratic Dirichlet
space to the $p$-Dirichlet scale, $1<p<\infty$.

Throughout the paper, $\Gamma$ is a rectifiable Jordan curve.  Its
critical Besov space $B_p(\Gamma)=B_{p,p}^{1/p}(\Gamma)$ consists, modulo
additive constants, of all measurable functions $u:\Gamma\to\mathbb C$ such
that
$$\|u\|_{B_{p}(\Gamma)}^p = \frac{1}{4\pi^2} \iint_{\Gamma\times \Gamma}  \frac{|u(z_1) - u(z_2)|^p}{|z_1 - z_2|^2} |dz_1||dz_2|<\infty.$$
 
By the Jordan curve theorem, $\Gamma$ separates the plane into the bounded
component $\Omega_i$ and the unbounded component $\Omega_e$.  For $p>1$, let
$\mathcal D_p(\Omega_i)$ be the space of harmonic functions $U_i$ on
$\Omega_i$ with finite $p$-Dirichlet energy
\begin{equation}
    D_p(U_i) = \frac{1}{2\pi}\iint_{\Omega_i}|\nabla U_i(w)|^p \lambda_{\Omega_i}^{2-p}(w) dA(w) < \infty.
\end{equation} 
Here, $dA$ is planar area measure and
\[
\lambda_{\Omega_i}(\varphi(z))
=\frac{1}{(1-|z|^2)|\varphi'(z)|},
\]
where $\varphi:\mathbb D_i\to\Omega_i$ is conformal from the unit disk $\mathbb D_i$ onto $\Omega_i$.  For a
complex-valued function $U$ and $w=u+iv$, we use the Hilbert--Schmidt norm
\begin{equation}\label{gradient-complex}
|\nabla U|=\bigl(|U_u|^2+|U_v|^2\bigr)^{1/2}
=\bigl(2|U_w|^2+2|U_{\bar w}|^2\bigr)^{1/2}.
\end{equation}
In particular, $|\nabla U|=|U_w|+|U_{\bar w}|$ when $U$ is
real-valued, whereas $|\nabla U|=\sqrt2\,|U'|$ when $U$ is holomorphic.

The $p$-Dirichlet energy is invariant under conformal and
anticonformal bijections.  Indeed, let $\Phi:\Omega_1\to\Omega_2$ be such a
bijection and define $\sigma_\Phi=|\Phi_z|$ in the conformal case and
$\sigma_\Phi=|\Phi_{\bar z}|$ in the anticonformal case.  Under the change of
variables $w=\Phi(z)$,
$$
dA(w) = \sigma_\Phi(z)^2 dA(z), \qquad \lambda_{\Omega_1}(z) = \lambda_{\Omega_2}(\Phi(z))\sigma_{\Phi}(z).
$$
For a harmonic function $U$ on $\Omega_2$, the chain rule gives 
 $$
 |\nabla(U\circ\Phi)(z)| = \sigma_\Phi(z)|\nabla U(\Phi(z))|.
 $$
Combining these identities and applying the change of variables gives
$D_p(U\circ\Phi)=D_p(U)$.

A function $U_i\in\mathcal D_p(\Omega_i)$ has conformally
nontangential boundary values $u$ at arclength-almost every point of $\Gamma$.
Moreover, it is recovered from $u$ by the Poisson formula
$$
U_i = P(u\circ\varphi)\circ\varphi^{-1},
$$
where $P$ is the classical Poisson operator on $\mathbb D_i$.  We call
$D_p(U_i)$ the interior $p$-Dirichlet energy of $u$ and write it as
$\Vert u\Vert_i^p$.

For completeness, Peller's result (\cite{pel}) asserts that  $\mathcal{D}_p(\mathbb{D}_i)\subset \text{VMO}(\mathbb{D}_i)$ for  $p >1$. Since the latter space is  included in the Hardy space of harmonic functions $h^2(\mathbb{D}_i)$, $U_i\circ\varphi \; (\in \mathcal D_p(\mathbb D_i))$ has angular limits almost everywhere on the unit circle $\mathbb S$.   The F.~and M.~Riesz
theorem for rectifiable Jordan curves states that harmonic measure and arclength measure are mutually absolutely continuous. It therefore transfers the almost-everywhere
statement between $\mathbb S$ and $\Gamma$.

To define $\mathcal D_p(\Omega_e)$, regard $\infty$ as an interior
point.  After translating and scaling if necessary, assume that
$0\in\Omega_i \subset\mathbb D_i$ and set $\iota(z)=1/\bar z$.  Then $\iota$ maps $\Omega_e$
onto a bounded Jordan domain $\widetilde\Omega_i$ and sends $\infty$ to $0$.
We call $U_e$ harmonic on $\Omega_e\cup\{\infty\}$ when
$\widetilde U_i:=U_e\circ\iota$ extends harmonically across $0$, and we set
$U_e(\infty)=\widetilde U_i(0)$.  Thus
\[
\mathcal D_p(\Omega_e)=
\{\widetilde U_i\circ\iota:\widetilde U_i\in
\mathcal D_p(\widetilde\Omega_i)\}.
\]
Its energy is
\[
D_p(U_e)=\frac{1}{2\pi}\iint_{\Omega_e}
|\nabla U_e(w)|^p\lambda_{\Omega_e}^{2-p}(w)\,dA(w)
=D_p(\widetilde U_i),
\]
where the point $\infty$ is omitted from the Euclidean integral and the
last equality follows from anticonformal invariance.

If $\widetilde\varphi:\mathbb D_i\to\widetilde\Omega_i$ is
conformal, then
\[
U_e=P(u\circ\iota\circ\widetilde\varphi)
\circ\widetilde\varphi^{-1}\circ\iota.
\]
We call $D_p(U_e)$ the exterior $p$-Dirichlet energy of $u$ and denote
it by $\Vert u\Vert_e^p$.

For $p=2$ and $\Gamma=\mathbb S$, Douglas proved that the boundary
space of $\mathcal D_2(\mathbb D_i)$ is $B_2(\mathbb S)$ and that
\begin{equation*}
    \frac{1}{2\pi}\iint_{\mathbb D_i}|\nabla U_i(w)|^2\,dA(w)
    =   \frac{1}{4\pi^2} \iint_{\mathbb S\times \mathbb S}  \frac{|u(z_1) - u(z_2)|^2}{|z_1 - z_2|^2} |dz_1||dz_2|.
\end{equation*}
This is the Douglas formula, introduced in Douglas's solution of the
Plateau problem \cite{dou}.  It also gives
$$
\frac{1}{4\pi^2} \iint_{\mathbb S\times \mathbb S}
\frac{|u(z_1) - u(z_2)|^2}{|z_1 - z_2|^2}\,|dz_1|\,|dz_2|
=\frac{1}{2\pi}\iint_{\mathbb D_e}|\nabla U_e(w)|^2\,dA(w).
$$
Thus the same boundary seminorm has three equal representations:
$$
\Vert u \Vert_i = \Vert u \Vert_{B_2(\mathbb S)} = \Vert u \Vert_e.
$$

In \cite{WZ} we identified the rectifiable Jordan curves for which
these three seminorms are equivalent:
\[
\Vert \cdot\Vert_i\simeq\Vert \cdot\Vert_{B_2(\Gamma)}
\simeq\Vert \cdot\Vert_e
\]
holds if and only if $\Gamma$ is a chord-arc curve.  Here and below,
$A\simeq B$ means that $C^{-1}A\le B\le CA$ for a constant independent of
the function under consideration.

Recall that a rectifiable Jordan curve is chord-arc (or
$K$-chord-arc) if the shorter of the two subarcs joining any
$z_1,z_2\in\Gamma$ has length at most $K|z_1-z_2|$.

For every $p>1$, it is known that
\[
\Vert\cdot\Vert_i\simeq\Vert\cdot\Vert_e
\simeq\Vert\cdot\Vert_{B_p(\Gamma)}
\]
on a chord-arc curve; see \cite{LS}; that is,  there
is a constant $C>0$ for which the following inequalities hold:
 \begin{itemize}
\item[(a)]$\|\cdot\|_{B_{p}(\Gamma)}\le C \|\cdot\|_{i}$ and $\|\cdot\|_{B_{p}(\Gamma)}\le C \|\cdot\|_{e}$;
\item[(b)]$\|\cdot\|_{i}\le C\|\cdot\|_{B_{p}(\Gamma)}$ and $\|\cdot\|_{e}\le C\|\cdot\|_{B_{p}(\Gamma)}$.
\end{itemize}
More precisely, the interior and
exterior energies are equivalent if and only if the curve $\Gamma$ is a quasicircle, by the critical
composition theorem \cite[Theorem~12]{Bou} and conformal and anticonformal invariances of $p$-Dirichlet energies. Thus, if $\Gamma$ is a rectifiable quasicircle then one inequality in \rm{(a)} (also in \rm{(b)}) implies the other with a different constant $C$. 

We investigate the converses to these estimates.
\begin{theorem}[Main Theorem] \label{rec} 
Let $\Gamma$ be a rectifiable quasicircle.  Then the following statements
hold.
\begin{itemize}
    \item[{\rm(1)}] If \rm{(a)} holds for some $p>1$, then $\Gamma$ is a chord-arc curve.
    \item[{\rm(2)}] There exists a non-chord-arc rectifiable quasicircle $\Gamma$ for which \rm{(b)} holds for every $p>1$.
\end{itemize}
\end{theorem}
\begin{remark}
    For a harmonic function $U$ on one side  $\Omega = \{\Omega_i, \Omega_e\}$,  define the weighted energy for $1<p<\infty$, $0<s<1$, as
$$
\mathcal E_{s,p;\Omega}(U) :=  \left(\iint_\Omega |\nabla U(z)|^p\lambda_\Omega(z)^{1-(1-s)p}dA(z)\right)^{1/p}.
$$
The $p$-Dirichlet spaces $\mathcal D_p(\Omega)$ are a special case of  Besov-Sobolev spaces
\[  
\mathcal B_{p,p}^s(\Omega)  
:=\{U:\ U\text{ harmonic in }\Omega,             
\ \mathcal E_{s,p;\Omega}(U)<\infty\},
\]
understood modulo constants.
Indeed, the $p$-Dirichlet spaces correspond to the case $s=1/p$ and these are precisely the conformally invariant members of the
Besov--Sobolev scale.  Whenever
the corresponding trace estimate holds, the trace operator is first defined
on smooth functions by restriction to $\Gamma$ and then extended to
$\mathcal B_{p,p}^s(\Omega)$ by completion. If $\Gamma=\partial\Omega$ is a rectifiable curve we equip boundary measurable functions $u:\Gamma\to\mathbb C$ with the
homogeneous Douglas-type Besov seminorm 
$$
\Vert u\Vert_{B_{p,p}^s(\Gamma)} := \left(\iint_{\Gamma\times\Gamma}\frac{|u(z)-u(\zeta)|^p}{|z-\zeta|^{1+sp}}|dz||d\zeta|\right)^{1/p}<\infty.
$$
and denote by $B_{p,p}^s(\Gamma)$ the space of measurable functions for which
this seminorm is finite, understood modulo constants.

We announce here that 
suppose $\Gamma$ is a rectifiable quasicircle then if $\Gamma$ is chord-arc then 
\begin{equation}    
\|\operatorname{Tr}_\Gamma U\|_{B^s_{p,p}(\Gamma)} \lesssim \mathcal E_{s,p;\Omega_i}(U), \qquad \|\operatorname{Tr}_\Gamma U\|_{B^s_{p,p}(\Gamma)} \lesssim \mathcal E_{s,p;\Omega_e}(U),
\end{equation}
and conversely, the validity of either one of these inequalities implies $\Gamma$ is chord-arc. This, in particular, shows that the converse assertion in part~{\rm(1)} of Main Theorem
is not a consequence of the conformal invariance peculiar to the critical
space $\mathcal D_p(\Omega)$. Rather, it is a geometric property shared by
the entire Besov--Sobolev scale considered above. A detailed proof,
together with a careful treatment of the more precise relationship between $B_{p,p}^s(\Gamma)$ and $\text{Tr}_\Gamma\mathcal B_{p,p}^s(\Omega)$, will appear in a
forthcoming paper.
\end{remark}


The paper is organized as follows. 
Section~\ref{sec:Ar} gives a characterization of Ahlfors-regular curves.  
Section~\ref{sec:bigp} uses this characterization to prove part~\rm{(1)} for $2\le p<\infty$, while Section~\ref{sec:smallp} treats $1<p<2$.  Section~\ref{sec:ex} constructs 
the example in part~\rm{(2)}.

\section{Ahlfors-regular curves}\label{sec:Ar}
 We begin with the definition of Ahlfors regularity.
\begin{definition} 
Let $\Gamma$ be a rectifiable curve in the plane.  We call $\Gamma$
$M$-regular if, for every $z\in\mathbb C$ and $r>0$,
$$\mathrm{length}(\Gamma\cap B(z,r))\leq Mr,$$
where $B(z, r)$ stands for the open disk centered at $z$ and of radius $r$. 
\end{definition}

\begin{theorem}\label{imp}
    Let $\Gamma$ be a rectifiable curve in the plane.  The following
    statements are equivalent.
\begin{itemize}
    \item[{\rm(i)}] $\Gamma$ is an $M$-regular curve; 
    \item[{\rm(ii)}] There exists $C>0$ such that for every $w\notin \Gamma$,
$$\int_{\Gamma}\frac{|dz|}{|z - w|^2}  \leq \frac{C}{d(w,\Gamma)},$$
where $d(w,\Gamma)$ denotes the distance from $w$ to $\Gamma$.
\end{itemize}
\end{theorem}
\noindent A proof of $\rm(ii)\Rightarrow\rm(i)$ was given in
\cite{zin}, where the result is attributed to Y.~Meyer; see also \cite{WZ}.
We give a different proof, based on porosity and inspired by \cite{bru}.
The elementary implication $\rm(i)\Rightarrow\rm(ii)$ is proved in
\cite{WZ,bru}.

\begin{definition} A compact set $K\subset\mathbb C$ is porous if there
exist constants $0<c<1$ and $r_0>0$ such that, for every $z\in K$ and every
$0<r\le r_0$, the disk $B(z,r)$ contains a disk $B(z',cr)$ disjoint from $K$.
\end{definition}
\begin{lemma}\label{poro} 
If $\Gamma$ satisfies \rm(ii) and is porous, then it
satisfies \rm(i).
\end{lemma}
\begin{proof}
It is enough to consider disks \(B(z,r)\) with \(z\in\Gamma\) and
\(0<r\le r_0\); larger disks are handled using the finite length of
\(\Gamma\). Let
\[
  B(z',cr)\subset B(z,r)\setminus\Gamma
\]
be the porous hole. For every \(w\in\Gamma\cap B(z,r)\), we have
\(|w-z'|\le2r\), while \(d(z',\Gamma)\ge cr\). Hence
\[
\begin{aligned}
 {\rm length}(\Gamma\cap B(z,r)) \le 4r^2
   \int_{\Gamma\cap B(z,r)}
   \frac{|dw|}{|w-z'|^2} \le \frac{4Cr^2}{d(z',\Gamma)}
 \le \frac{4C}{c}\,r.
\end{aligned}
\]
This proves (i).
\end{proof}

\begin{proof}[Proof of {\rm (ii) $\Rightarrow$ (i)}]
For convenience, we use the \(\ell^\infty\)-norm
\[
  \|x+iy\|_\infty=\max\{|x|,|y|\},
\]
so that \(B_\infty(z_0,r)\) is the square centered at \(z_0\) with side
length \(2r\). Since the Euclidean and \(\ell^\infty\)-norms are
equivalent, this change affects only the constants.

By Lemma~\ref{poro}, it is enough to prove that condition (ii)
implies that \(\Gamma\) is porous. Consider the closed square
\(\overline{B_\infty(z_0,r)}\), and choose
\(\zeta\in\overline{B_\infty(z_0,r)}\) such that
\[
 d(\zeta,\Gamma)
 =\max\bigl\{d(z,\Gamma):
 z\in\overline{B_\infty(z_0,r)}\bigr\}.
\]
There is a square contained in \(B_\infty(z_0,r)\), having \(\zeta\) as
a vertex and side length \(d(\zeta,\Gamma)\), whose interior is disjoint
from \(\Gamma\). Let \(z_0'\) be the center of this square. Thus the
square is
\[
 B_\infty\left(z_0',\frac{d(\zeta,\Gamma)}{2}\right).
\]
If
\[
 \delta=d(z_0',\Gamma),
\]
then
\begin{equation}\label{eq:21}
 \frac12d(\zeta,\Gamma)
 \le\delta
 \le\frac32d(\zeta,\Gamma).
\end{equation}

By construction, every closed subsquare of
\(\overline{B_\infty(z_0,r)}\) with side length \(2\delta\) meets
\(\Gamma\). To prove porosity, it remains to show that
\(d(\zeta,\Gamma)\) is bounded below by a constant multiple of \(r\).
In view of \eqref{eq:21}, it is enough to prove the corresponding lower
bound for \(\delta\).

Condition (ii), applied at \(z_0'\), gives
\[
  \int_\Gamma\frac{|dz|}{|z-z_0'|^2}
  \le\frac{C}{\delta}.
\]
We next obtain a lower bound for this integral. By translating and
rescaling, we may assume that the square under consideration is
\(B_\infty(0,1)\). We may also assume that
\(\diam(\Gamma)\ge3\).

The square \(B_\infty(0,1)\) contains a square \(\mathcal C\) of side
length \(1\) having \(z_0'\) as a vertex. We define the \emph{heart} of
a square \(B_\infty(\xi,\rho)\) to be the concentric square
\(B_\infty(\xi,\rho/3)\). If \(\Gamma\) meets the heart of
\(B_\infty(\xi,\rho)\), then
\[
 {\rm length}\bigl(\Gamma\cap B_\infty(\xi,\rho)\bigr)
 \ge \frac{2\rho}{3}.
\]

Divide \(\mathcal C\) into \(N^2\) subsquares of side length \(1/N\),
where \(N\) is chosen so that
\begin{equation}\label{eq:22}
 \frac{1}{3N}\ge2\delta\ge\frac{1}{3(N+1)}.
\end{equation}
The right-hand inequality in \eqref{eq:22} shows that the Euclidean
diagonal \(\sqrt2/N\) of each subsquare is less than \(20\delta\),
provided that \(\delta<1/12\). If \(\delta\ge1/12\), the desired lower
bound has already been obtained. The left-hand inequality in
\eqref{eq:22} shows that the side length of each subsquare is at least
\(6\delta\).

Let \(\mathcal C_j\) be one of these subsquares whose distance from
\(z_0'\) is at least \(6\delta\). The preceding observations imply that
\(\Gamma\) meets the heart of \(\mathcal C_j\); consequently,
\({\rm length}(\Gamma\cap\mathcal C_j)\ge1/(4N)\). It follows that
\[
\begin{aligned}
 \int_{\Gamma\cap\mathcal C_j}
 \frac{|dz|}{|z-z_0'|^2}
 &\ge
 \frac{1}
 {4N\bigl(\sup_{z\in\mathcal C_j}|z-z_0'|\bigr)^2} \\
 &\ge
 \frac{N}{10}
 \iint_{\mathcal C_j}
 \frac{dA(z)}{|z-z_0'|^2}.
\end{aligned}
\]
Summing over all such subsquares \(\mathcal C_j\), we obtain
\[
\begin{aligned}
 \sum_j\int_{\Gamma\cap\mathcal C_j}
 \frac{|dz|}{|z-z_0'|^2}
 &\ge
 \frac{N}{10}\,\frac{\pi}{2}
 \int_{20\delta}^{1}\frac{dr}{r} \\
 &\ge
 \frac{N}{10}\log\frac{1}{20\delta} \\
 &\ge
 \frac{1}{1000\delta}
 \log\frac{1}{20\delta}.
\end{aligned}
\]
Combining this lower bound with the preceding upper estimate gives
\[
 \frac{1}{1000\delta}
 \log\frac{1}{20\delta}
 \le\frac{C}{\delta}.
\]
Therefore
\[
 \delta\ge\frac1{20}e^{-1000C}.
\]
Here \(C\) is the constant in condition (ii). Together with
\eqref{eq:21}, this gives the required porous hole and completes the
proof.
\end{proof}

\section{Part (1) of Theorem \ref{rec}: case
\texorpdfstring{$p\geq 2$}{p >= 2}}\label{sec:bigp}
Fix $w\in\Omega_e$ and write $d=d(w,\Gamma)$.  If
$d\ge\frac14\diam(\Gamma)$, then
\begin{equation}\label{ddd}
    \int_{\Gamma}\frac{|dz|}{|z - w|^2} \leq \frac{C'}{d}
\end{equation}
where $C'=4\operatorname{length}(\Gamma)/\diam(\Gamma)$.  We now
assume that $d<\frac14\diam(\Gamma)$ and prove the same estimate with a
different constant.

Let $U_i(z)=(z-w)^{-1}$ on $\Omega_i$ and let $u=U_i|_\Gamma$.
Koebe's estimate
\[
\lambda_{\Omega_i}(z)\le \delta(z)^{-1}
\le4\lambda_{\Omega_i}(z),\qquad
\delta(z)=d(z,\Gamma),
\]
together with \eqref{gradient-complex}, gives
 \begin{align*}
   \Vert u \Vert_i^p
   &= \frac{2^{p/2}}{2\pi}\iint_{\Omega_i}|U_i'(z)|^p
   \lambda_{\Omega_i}^{2-p}(z)\,dA(z)\\
   &\leq \frac{2^{p/2}4^{p-2}}{2\pi}
   \iint_{\Omega_i}\frac{\delta(z)^{p-2}}{|z-w|^{2p}}\,dA(z)\\
   &\leq \frac{2^{p/2}4^{p-2}}{2\pi}
   \iint_{|z-w|\ge d}|z-w|^{-p-2}\,dA(z)\\
   &=\frac{2^{p/2}4^{p-2}}{p}\,d^{-p}.
\end{align*}
Here, $p\ge2$ is used when replacing
$\lambda_{\Omega_i}^{2-p}$ by $4^{p-2}\delta^{p-2}$ and when using
$\delta(z)\le|z-w|$.  Choose $z_0,z_1\in\Gamma$ so that
$|w-z_0|=d$ and $|w-z_1|=D:=\max_{z\in\Gamma}|w-z|$.  Then
$D\ge\frac12\diam(\Gamma)$.  Under the map $z\mapsto(z-w)^{-1}$, set
$\eta_j=(z_j-w)^{-1},\;j=0,1$ and denote the image of $\Gamma$ by $\Gamma'$.
We have
 $$
 |\eta_0 - \eta_1| \geq \frac{1}{d} - \frac{1}{D} \geq \frac{1}{d} - \frac{2}{\text{diam}(\Gamma)} \geq \frac{1}{d} - \frac{1}{2d} = \frac{1}{2d}. 
 $$For every $\eta\in\Gamma'$, the triangle inequality shows that at least
one of $|\eta-\eta_0|$ and $|\eta-\eta_1|$ is at least $1/(4d)$.
Following $\Gamma'$ from $\eta$ to that point, the distance from $\eta$
crosses the interval $[1/(8d),1/(4d)]$; the corresponding subarc has length
at least $1/(8d)$.  This elementary point, illustrated in Figure~\ref{fig:inversion},
justifies the estimate below.
\begin{align*}
    4\pi^2\Vert u \Vert_{B_p(\Gamma)}^p
    &=\int_{\Gamma}\int_{\Gamma}
    \frac{\left|(z-w)^{-1}-(\zeta-w)^{-1}\right|^p}{|z-\zeta|^2}
    \,|dz|\,|d\zeta|\\
    &=\int_{\Gamma'}\int_{\Gamma'}|\eta-\xi|^{p-2}
    \,|d\eta|\,|d\xi|\\
    &\ge \int_{\Gamma'}\!|d\eta|
    \int_{\{\xi\in\Gamma':\,1/(8d)<|\eta-\xi|<1/(4d)\}}
    |\eta-\xi|^{p-2}\,|d\xi|\\
    &\ge 8^{1-p}\operatorname{length}(\Gamma')d^{1-p}.
\end{align*}
The last inequality again uses $p\ge2$.  Notice also that
\begin{equation}\label{length-image}
\operatorname{length}(\Gamma')
=\int_\Gamma\frac{|dz|}{|z-w|^2}.
\end{equation}

Let $C$ be the constant in the interior inequality in \rm(a).
Combining the upper and lower bounds with
$\Vert u\Vert_{B_p(\Gamma)}\le C\Vert u\Vert_i$ and using
\eqref{length-image}, we obtain
\begin{equation}\label{dd}
    \int_{\Gamma}\frac{|dz|}{|z - w|^2} \leq \frac{C'}{d}
\end{equation}
where one may take
$C'=4\pi^2C^p8^{p-1}2^{p/2}4^{p-2}/p$.
\begin{figure}[htp]
\centering
\begin{tikzpicture}[
    scale=1.05,
    >=Latex,
    point/.style={
        circle,
        fill=black,
        inner sep=1.7pt
    },
    every node/.style={font=\small}
]
    \def\R{2.75}

    \fill[white!3] (0,0) circle (\R);
    \draw[black!55,line width=.9pt]
        (0,0) circle (\R);

    \node[black] at (-1.55,2.30)
        {$|\eta|=1/d$};

    \coordinate (e0) at (2.75,0);
    \coordinate (e1) at (.42,-.10);
    \coordinate (e)  at (.68,.62);

    \filldraw[
        fill=white!11,
        draw=black,
        line width=1.35pt
    ]
        (e0)
        .. controls (2.55,-.48) and (2.08,-1.05)
            .. (1.30,-1.12)
        .. controls (.78,-1.18) and (.44,-.65)
            .. (e1)
        .. controls (.38,.15) and (.45,.40)
            .. (e)
        .. controls (.96,1.05) and (1.48,1.18)
            .. (2.05,.85)
        .. controls (2.52,.62) and (2.70,.30)
            .. cycle;

    \node[
        black,
        fill=white,
        inner sep=1pt
    ] at (1.55,1.00) {$\Gamma'$};

    \fill[white] (0,0) circle (1.3pt);
    \node[black,anchor=north east]
        at (-.05,-.04) {$0$};

    \node[
        point,
        label={[black]right:$\eta_0$}
    ] at (e0) {};

    \node[
        point,
        label={[black]below left:$\eta_1$}
    ] at (e1) {};

    \node[
        point,
        label={[black]above left:$\eta$}
    ] at (e) {};

    \draw[
        black!75,
        dashed,
        line width=.85pt
    ]
        (e0)--(e1)
        node[
            midway,
            below=2pt,
            inner sep=1.5pt
        ]
        {};

    \draw[black!45,densely dotted]
        (e)--(e0);

    \draw[black!45,densely dotted]
        (e)--(e1);

    \node[
        align=center,
        black,
        inner sep=2pt
    ] at (-1.35,-1.82)
        {$\max\{|\eta-\eta_0|,|\eta-\eta_1|\}
          \ge 1/(4d)$};

    \draw[
        -{Latex[length=2mm]},
        black!65
    ]
        (-.25,-1.65)
        to[bend left=13]
        (.34,-.30);
\end{tikzpicture}

\caption{
The inverted curve $\Gamma'$ lies in
$\{\eta:|\eta|\le 1/d\}$.
The origin belongs to the unbounded component of
$\mathbb C\setminus\Gamma'$, and $\Gamma'$ contains two points $\eta_0$, $\eta_1$
separated by at least $1/(2d)$.
}
\label{fig:inversion}
\end{figure}
Combining \eqref{ddd} and \eqref{dd}, we obtain a constant $C'$ such
that
\begin{equation}\label{dddd}
    \int_{\Gamma}\frac{|dz|}{|z - w|^2} \leq \frac{C'}{d(w, \Gamma)}
\end{equation}
for every $w\in\Omega_e$, provided that the interior inequality in
\rm(a) holds.

For $w\in\Omega_i$, apply the exterior inequality in \rm(a) to
$U_e(z)=(z-w)^{-1}$ on $\Omega_e$.  The same argument proves
\eqref{dddd}.  Theorem~\ref{imp} now implies that $\Gamma$ is Ahlfors-regular.  An Ahlfors-regular quasicircle is chord-arc, so part~\rm(1) of
Theorem~\ref{rec} follows for $p\ge2$.

\section{Part (1) of Theorem \ref{rec}: case
\texorpdfstring{$1<p<2$}{1 < p < 2}}\label{sec:smallp}

Throughout this section, $\mathbb D = \mathbb D_i$, $\Omega=\Omega_i$ and $\delta(x)=\dist(x,\Gamma)$ for $x\in\Omega$.

We  use the bounded-turning characterization of quasicircles: a Jordan curve $\Gamma$ is
a $K$-quasicircle if, for every $a,b\in\Gamma$, one of the two subarcs
joining $a$ to $b$ has diameter at most $K|a-b|$.

We  prove part~\rm(1) of Theorem~\ref{rec} for $1<p<2$ using the
interior inequality in \rm(a), bounded turning, and the almost-Dirichlet
principle in Lemma~\ref{almost1p}.  It suffices to prove that there is a
constant $M$, independent of $x\in\Gamma$ and $0<r<\diam(\Gamma)$, such
that
\begin{equation}\label{xy}
    {\rm length}(\Gamma\cap B(x, r))\leq Mr.
\end{equation} 
Indeed, this implies Ahlfors regularity.  If
$B(y,r)\cap\Gamma\ne\varnothing$, choose $x\in B(y,r)\cap\Gamma$; then
$\Gamma\cap B(y,r)\subset\Gamma\cap B(x,2r)$, and \eqref{xy} gives
$$
{\rm length}(\Gamma\cap B(y, r))\leq 2Mr.
$$
The case $B(y,r)\cap\Gamma=\varnothing$ is trivial.  Thus $\Gamma$ is
Ahlfors-regular and, being a quasicircle, is chord-arc.

For $x\in\Gamma$ and $r>0$, write
$L_x(r)=\operatorname{length}(\Gamma\cap B(x,r))$.
\begin{lemma}\label{2arc}
Let $\Gamma$ be a rectifiable $K$-quasicircle of length $\ell$.  For
every $x\in\Gamma$ and $R>0$, the set $\Gamma\cap B(x,R)$ is contained in
the union of two 
subarcs $\alpha_+$ and $\alpha_-$ issuing from $x$, such that
\[
  \alpha_+\cup\alpha_-\subset B(x,2KR)
\]
and
\begin{equation}\label{alph}
  {\rm length}(\alpha_+)+{\rm length}(\alpha_-)
  \le2L_x(2KR).
\end{equation}
\end{lemma}

\begin{proof}
Choose an arclength parametrization
\[
  \gamma:[0,\ell]\longrightarrow\Gamma,
  \qquad \gamma(0)=\gamma(\ell)=x,
\]
which is one-to-one on $[0,\ell)$. For $s\in(0,\ell)$, the two arcs joining
$x$ and $\gamma(s)$ are $ \gamma([0,s])$ and $\gamma([s,\ell]).$  
If $\gamma(s)\in B(x,R)$, then at least one of these two
arcs has diameter at most $K|x-\gamma(s)|<KR$, so this arc is contained in
$B(x,KR)$.

The parameters $s$ for which $\diam(\gamma([0,s]))<KR$ form an
initial interval; its closure determines $\alpha_+$.  Define $\alpha_-$
similarly from the opposite orientation.  Both arcs lie in $B(x,2KR)$, and
the preceding paragraph shows that their union contains
$\Gamma\cap B(x,R)$.  Since each $\alpha_\pm$ is contained in
$\Gamma\cap B(x,2KR)$,
\[
  {\rm length}(\alpha_\pm)\le L_x(2KR).
\]
Adding the two estimates proves \eqref{alph}.
\end{proof}

We use the following tubular-neighborhood estimate
\cite[Lemma~4.2]{MT}: if $\alpha$ is a rectifiable arc, then
\begin{equation}\label{elem}
  \bigl|\{w:\dist(w,\alpha)<t\}\bigr|
  \le C\bigl(t\ell(\alpha)+t^2\bigr)
\end{equation}
where $C$ is universal.  This can also be seen directly:
choose at most $\lfloor\ell(\alpha)/t\rfloor+2$ points with successive
arclength separation at most $t$.  Disks of radius $2t$ centered at these
points cover the $t$-neighborhood and yield \eqref{elem}.

\begin{lemma}\label{8Kr}
Let $1<p<2$ and let $\Gamma$ be a rectifiable $K$-quasicircle.  Then,
for every $x\in\Gamma$ and $r>0$,
\begin{equation}\label{Lx}
  \int_{\Omega\cap B(x,2r)}\delta(w)^{p-2}\,dA(w)
  \leq C_p\bigl(r^{p-1}L_x(8K r)+r^p\bigr). 
\end{equation}
\end{lemma}

\begin{proof}
For $x\in\Gamma$ and $0<t\le2r$, define
\[
  E_t=\{w\in\Omega\cap B(x,2r):\delta(w)<t\}.
\]
If $w\in E_t$, choose $\xi\in\Gamma$ with
$\delta(w)=|\xi-w|$.  Then
\[
  |\xi-x|
  \le|\xi-w|+|w-x|
  =\delta(w)+|w-x|
  <t+2r\le4r.
\]
Lemma~\ref{2arc}, applied with $R=4r$, shows that all such points $\xi$ are
contained in two arcs $\alpha_+$ and $\alpha_-$ lying in
$B(x,8Kr)$ with
\[
  {\rm length}(\alpha_+)+ {\rm length}(\alpha_-)
  \le2L_x(8K r).
\]
Consequently, $E_t$ lies in the $t$-neighborhood of
$\alpha_+\cup\alpha_-$.  By \eqref{elem},
\begin{equation}\label{Et}
  |E_t|\le C\bigl(tL_x(8K r)+t^2\bigr),
  \qquad 0<t\le2r.
\end{equation}

Every $w\in B(x,2r)$ satisfies $\delta(w)\le|w-x|<2r$.  Put
$t_k=2^{1-k}r$ and decompose $\Omega\cap B(x,2r)$ into the layers
\[
  D_k=\{w:t_{k+1}\le\delta(w)<t_k\},
  \qquad k=0,1,2,\cdots.
\]
Since $D_k\subset E_{t_k}$ and, because $p-2<0$,
\[
\delta(w)^{p-2}\le t_{k+1}^{p-2}=2^{2-p}t_k^{p-2}
\quad\text{on }D_k,
\]
estimate \eqref{Et} gives
\begin{align*}
  \int_{\Omega\cap B(x,2r)}\delta(w)^{p-2}\,dA(w)
  &\le C_p\sum_{k=0}^\infty
  t_k^{p-2}\bigl(t_kL_x(8K r)+t_k^2\bigr)\\
  &=C_p\left(
  L_x(8K r)\sum_{k=0}^\infty t_k^{p-1}
  +\sum_{k=0}^\infty t_k^{p}
  \right).
\end{align*}
Both geometric series converge; their sums are bounded by constant
multiples of $r^{p-1}$ and $r^p$, respectively.  This proves \eqref{Lx}.
\end{proof}

For $p\ne2$, the harmonic extension need not minimize the weighted 
$p$-energy, but it satisfies an almost-Dirichlet principle.  Namely, if $U$
is the harmonic extension of $u$ and $V\in C(\overline\Omega)\cap
C^1(\Omega)$ is any extension of the same boundary data, then
$$
\iint_\Omega|\nabla U(z)|^p \delta(z)^{p-2}dA(z) \leq C \iint_\Omega|\nabla V(z)|^p \delta(z)^{p-2}dA(z)
$$
with a constant depending only on $p$.  The disk version
is the periodic counterpart of \cite[Theorem~1.11]{MR}.

\begin{lemma}\label{almost1p}
Let $\Omega$ be a bounded Jordan domain and $1<p<\infty$.  Then the
almost-Dirichlet principle holds for $\mathcal D_p(\Omega)$.
\end{lemma}
\begin{proof}
On $\mathbb D$, the principle reads
\begin{equation}    
\iint_{\mathbb D}|\nabla U(z)|^p (1-|z|)^{p-2}dA(z) \leq  C\iint_{\mathbb D}|\nabla V(z)|^p (1-|z|)^{p-2}dA(z) 
\end{equation}
where $U$ is harmonic and $V$ has the same boundary trace.  Let
$\varphi:\mathbb D\to\Omega$ be conformal.  Under $w=\varphi(z)$, the
weighted integral transforms as
\begin{align*}
&\int_{\mathbb D}|\nabla(F\circ\varphi)(z)|^p
(1-|z|)^{p-2}\,dA(z)\\
&\qquad=\int_\Omega |\nabla F(w)|^p
\bigl((1-|\varphi^{-1}(w)|)|\varphi'(\varphi^{-1}(w))|\bigr)^{p-2}
\,dA(w).
\end{align*}
Koebe's theorem gives
\[
(1-|z|)|\varphi'(z)|\simeq\delta(\varphi(z)),
\]
with absolute comparison constants.  Applying the disk principle to
$U\circ\varphi$ and $V\circ\varphi$, and then using this comparison, yields
\begin{equation}    
\iint_{\Omega}|\nabla U(w)|^p \delta(w)^{p-2}\,dA(w)
\leq C_p\iint_{\Omega}|\nabla V(w)|^p
\delta(w)^{p-2}\,dA(w).
\end{equation}
This proves the assertion.
\end{proof}

\medskip

We now prove part~\rm(1) of Theorem~\ref{rec} for $1<p<2$ in two
steps.

\textbf{Step 1: establish a recurrence for $L_x(r)/r$.}

Fix $x\in\Gamma$ and $r>0$. Choose a real-valued function
$\chi\in C_c^\infty(B(0,2))$ such that $\chi\equiv 1$ on $B(0,1)$, and define
\begin{equation}
  \begin{aligned}
  F_1(z)&=\chi\left(\frac{z-x}{r}\right)
  \frac{\operatorname{Re}(z-x)}{r},\\
  F_2(z)&=\chi\left(\frac{z-x}{r}\right)
  \frac{\operatorname{Im}(z-x)}{r}.
  \end{aligned}
\end{equation}
Then
\begin{equation}\label{supp}
  |\nabla F_j|\le\frac{C_\chi}{r},
  \qquad
  \operatorname{supp}(\nabla F_j)\subset B(x,2r),
  \qquad j=1,2.
\end{equation}
Let $u_j=F_j|_\Gamma$.  Since $\chi=1$ on $B(0,1)$, for
$z,\zeta\in \Gamma\cap B(x,r)$ we have
\[
  u_1(z)-u_1(\zeta)=\frac{\operatorname{Re}(z-\zeta)}{r},
  \qquad
  u_2(z)-u_2(\zeta)=\frac{\operatorname{Im}(z-\zeta)}{r}.
\]
Using
\[
  |\operatorname{Re}(z-\zeta)|^p
  +|\operatorname{Im}(z-\zeta)|^p
  \ge c_p|z-\zeta|^p
\]
and restricting both Besov integrals to
$(\Gamma\cap B(x,r))^2$, we obtain
\begin{equation}\label{sum}
 \lVert u_1\rVert_{B_p(\Gamma)}^p +  \lVert u_2\rVert_{B_p(\Gamma)}^p
  \ge c_pr^{-p}\int_{\Gamma\cap B(x,r)}\int_{\Gamma\cap B(x,r)} |z-\zeta|^{p-2}
  \,|dz|\,|d\zeta|.
\end{equation}
Since $|z-\zeta|\le2r$ and $p-2<0$,
\[
|z-\zeta|^{p-2}\ge(2r)^{p-2}.
\]
This is precisely where $p<2$ is used.  Thus \eqref{sum} implies
\begin{equation}\label{oneside}
 \lVert u_1\rVert_{B_p(\Gamma)}^p +  \lVert u_2\rVert_{B_p(\Gamma)}^p 
  \ge c_p\left(\frac{L_x(r)}{r}\right)^2.
\end{equation}
On the other hand, Lemma~\ref{almost1p}, \eqref{supp},
Lemma~\ref{8Kr}, and Koebe's comparison between the hyperbolic and distance
weights give
\begin{equation}\label{otherside}
    \lVert u_1\rVert_i^p + \lVert u_2\rVert_i^p  
\le Cr^{-p}\int_{\Omega\cap B(x,2r)}  
\delta(w)^{p-2}\,dA(w)   
\le C\left(1+\frac{L_x(8K r)}{r}\right).
\end{equation} 
By the interior inequality in \rm(a), with constant $C_0$,
\[
 \lVert u_1\rVert_{B_p(\Gamma)}^p + \lVert u_2\rVert_{B_p(\Gamma)}^p
  \le C_0^p (\lVert u_1\rVert_i^p + \lVert u_2\rVert_i^p).
\]
Combining this estimate with \eqref{oneside} and \eqref{otherside}
yields
\begin{equation}\label{recu}
  \left(\frac{L_x(r)}{r}\right)^2
  \le C_1\left(1+\frac{L_x(8K r)}{r}\right),
\end{equation}
where $C_1$ depends only on $p$, $C_0$, $K$, and the fixed cutoff $\chi$.

\textbf{Step 2: iterate to prove $L_x(r)\lesssim r$.}

Put
\[
  A_x(r)=\frac{L_x(r)}{r}.
\]
Since
\[
  \frac{L_x(8K r)}{r}=8K A_x(8K r),
\]
we may absorb $8K$ into the constant in \eqref{recu} and obtain
\begin{equation}
  A_x(r)^2\le C_2\bigl(1+A_x(8K r)\bigr).
\end{equation}
After increasing a constant $B\ge1$ if necessary,
\begin{equation}\label{root-recurrence}
  A_x(r)\le B\max\{1,A_x(8K r)^{1/2}\}.
\end{equation}

Let $D=\diam(\Gamma)$ and $L=\operatorname{length}(\Gamma)$.
If $r\ge D$, then $A_x(r)\le L/D$.  If $0<r<D$, let $N$ be the least
positive integer such that $(8K)^Nr\ge D$.  Applying
\eqref{root-recurrence} successively at the scales
$r,8Kr,\ldots,(8K)^{N-1}r$ gives
\begin{align}
  A_x(r)
  &\le B^{1+1/2+\cdots+1/2^{N-1}}
  \max\{1,A_x((8K)^Nr)^{1/2^N}\}\notag\\
  &\le B^2\max\left\{1,\frac{L}{D}\right\}.
\end{align}
Indeed, the geometric sum in the exponent is less than $2$, and
\[
  A_x((8K)^Nr)
  =\frac{L_x((8K)^Nr)}{(8K)^Nr}
  \le\frac{L}{D}.
\]
Thus $L_x(r)\le Mr$ with $M$ independent of $x$ and $r$.  This
completes the proof for $1<p<2$.

\section{Part \texorpdfstring{$(2)$}{(2)} of Theorem \ref{rec}}\label{sec:ex}
We now construct the example in part~\rm(2) of Theorem~\ref{rec}.
Let $\varphi:\mathbb D_i\to\Omega_i$ be a Riemann map.  The classical
rectifiability criterion states that $\Gamma=\partial\Omega_i$ is rectifiable
if and only if $\varphi'\in H^1(\mathbb D_i)$; see \cite[Chapter~3]{pom}.
\begin{definition} 
A rectifiable Jordan domain $\Omega_i$ is a Smirnov domain if
$\varphi'$ is outer in $H^1(\mathbb D_i)$; equivalently,
\begin{equation}\label{outer}
    \log|\varphi'(z)|
    =\frac{1}{2\pi}\int_{\mathbb S}
    \frac{1-|z|^2}{|\zeta-z|^2}\log|\varphi'(\zeta)|\,|d\zeta|,
    \qquad z\in\mathbb D_i.
\end{equation}
\end{definition}
Lavrentiev proved that chord-arc domains are Smirnov domains
\cite{lav}; more generally, a Jordan domain with Ahlfors-regular boundary is
Smirnov \cite{zinRegular}.  On the other hand, the constructions in \cite{dss,kah}
give a noncircular quasidisk $\Omega_i$ with rectifiable boundary for which
the derivative of a Riemann map is a nonconstant singular inner function:
\begin{equation}\label{ex}
    \vert\varphi'(\zeta)\vert = 1 \quad \text{\rm for almost\;all }\, \zeta \in \mathbb S, \quad (0 <) \vert\varphi'(z)\vert < 1 \quad \text{\rm for}\, z \in \mathbb D_i. 
\end{equation}
Consequently, harmonic measure at $\varphi(0)$ equals normalized
arclength measure.  Since a function that is both inner and outer is constant
up to a unimodular factor, \eqref{ex} and the noncircularity of $\Omega_i$
show that \eqref{outer} fails.  In particular, $\Gamma$ cannot be chord-arc.

Because $|\varphi'|=1$ almost everywhere on $\mathbb S$, the boundary
map $t\mapsto\varphi(e^{it})$ is an arclength parametrization and
$\operatorname{length}(\Gamma)=2\pi$.  Moreover, the shorter parameter arc
gives
\begin{equation}\label{upper-param}
|\varphi(e^{it})-\varphi(e^{is})|
\le \pi\left|\sin\frac{t-s}{2}\right|
=\frac{\pi}{2}|e^{it}-e^{is}|.
\end{equation}
It follows that
\begin{align*}
\Vert u \Vert_{B_p(\Gamma)}^p
&= \frac{1}{4\pi^2}\int_0^{2\pi}\int_0^{2\pi}
\frac{|u\circ\varphi(e^{it}) - u\circ\varphi(e^{is})|^p}
{|\varphi(e^{it})-\varphi(e^{is})|^2}\,dt\,ds\\
&\geq \frac{1}{\pi^4} \int_0^{2\pi}\int_0^{2\pi}
\frac{|u\circ\varphi(e^{it}) - u\circ\varphi(e^{is})|^p}
{|e^{it}-e^{is}|^2}\,dt\,ds\\
&= \frac{4}{\pi^2}\Vert u\circ\varphi \Vert_{B_p(\mathbb S)}^p.
\end{align*}
The unit-circle trace theorem (the periodic form of
\cite[Theorems~1.3 and~1.9]{MR}) asserts $\Vert u\circ\varphi\Vert_i \simeq \Vert u\circ\varphi\Vert_{B_p(\mathbb S)}$.  Conformal invariance then yields
\[
\Vert u\Vert_i\lesssim
\Vert u\circ\varphi\Vert_{B_p(\mathbb S)}
\lesssim\Vert u\Vert_{B_p(\Gamma)}.
\]
Finally, because $\Gamma$ is a quasicircle, the interior and exterior
$p$-Dirichlet energies are equivalent.  Hence
\[
\Vert u\Vert_e\lesssim\Vert u\Vert_i
\lesssim\Vert u\Vert_{B_p(\Gamma)}.
\]
Thus both inequalities in \rm(b) hold for every $p>1$, while $\Gamma$
is not chord-arc.  This proves part~\rm(2).

\paragraph*{Acknowledgments.}
The authors are grateful to Katsuhiko Matsuzaki for encouraging them
to complete the proof in the range $1<p<2$, presented in
Section~\ref{sec:smallp}, and for carefully reading the argument.
The recent preprint \cite{LSY} also announces a proof of the result in
this range. The proof presented here was obtained independently of
that work.

\begingroup
\small
\linespread{0.96}\selectfont

\endgroup

\end{document}